\documentclass[11pt]{article}

\usepackage{amsmath}
\usepackage{amsthm}
\usepackage{graphicx}
\usepackage{amsfonts}
\usepackage{amssymb}
\usepackage{psfrag}
\usepackage{wrapfig}
\usepackage{url}
\usepackage{amscd}
\usepackage{multicol}
\newcommand{\re}{\mathbb{R}}
\newcommand{\cpx}{\mathbb{C}}

\newcommand{\F}{\mathbb{F}}
\newcommand{\SF}{\mathbb{SF}}

\newcommand{\lmd}{\lambda}
\newcommand{\half}{\frac{1}{2}}

\newcommand{\af}{\alpha}

\newcommand{\reff}[1]{(\ref{#1})}

\newcommand{\prob}{\mbox{\rm Prob}\, }
\newcommand{\ex}{\mbox{\bf\sf E}}
\newcommand{\sign}{\mbox{\rm sign}}
\newcommand{\rank}{\mbox{\rm rank}}

\newcommand{\Var}{\mbox{\rm Var}}
\newcommand{\Tr}{\mbox{\rm Tr\,}}
\newcommand{\st}{\mbox{\rm s.t.}}
\newcommand{\II}{{\cal I}}
\newcommand{\DD}{{\cal D}}

\newcommand{\bdes}{\begin{description}}
\newcommand{\edes}{\end{description}}
\newcommand{\bal}{\begin{align}}
\newcommand{\eal}{\end{align}}
\newcommand{\bnum}{\begin{enumerate}}
\newcommand{\enum}{\end{enumerate}}
\newcommand{\bit}{\begin{itemize}}
\newcommand{\eit}{\end{itemize}}
\newcommand{\bea}{\begin{eqnarray}}
\newcommand{\eea}{\end{eqnarray}}
\newcommand{\be}{\begin{equation}}
\newcommand{\ee}{\end{equation}}
\newcommand{\baray}{\begin{array}}
\newcommand{\earay}{\end{array}}
\newcommand{\bsry}{\begin{subarray}}
\newcommand{\esry}{\end{subarray}}
\newcommand{\bca}{\begin{cases}}
\newcommand{\eca}{\end{cases}}
\newcommand{\bcen}{\begin{center}}
\newcommand{\ecen}{\end{center}}
\newcommand{\bbm}{\begin{bmatrix}}
\newcommand{\ebm}{\end{bmatrix}}
\newcommand{\bmx}{\begin{matrix}}
\newcommand{\emx}{\end{matrix}}
\newcommand{\bpm}{\begin{pmatrix}}
\newcommand{\epm}{\end{pmatrix}}
\newcommand{\btab}{\begin{tabular}}
\newcommand{\etab}{\end{tabular}}

\theoremstyle{plain}
\newtheorem{thm}{Theorem}[section]
\newtheorem{theorem}[thm]{Theorem}

\newtheorem{lemma}[thm]{Lemma}

\theoremstyle{definition}
\newtheorem{example}[thm]{Example}

\begin{document}

\title{ %On the Semidefinite Relaxation Bounds for Homogeneous Quadratic Optimization with Indefinite Matrix
Semidefnite Relaxation Bounds for Indefinite Homogeneous Quadratic Optimization
\author{
Simai He\footnote{Department of Systems Engineering and
Engineering Management, The Chinese University of Hong Kong,
Shatin, Hong Kong. Email:\,\texttt{smhe@se.cuhk.edu.hk}. }, \,
Zhi-Quan Luo\footnote{Department of Electrical and Computer
Engineering, University of Minnesota, 200 Union street SE,
Minneapolis, MN 55455. Email:\,\texttt{luozq@ece.umn.edu}. Research supported
in part by U.S.\
NSF grants DMS-0312416 and DMS-0610037.}, \,
Jiawang Nie\footnote{Institute of Mathematics and its
Applications, University of Minnesota, 207 Church street SE,
Minneapolis, MN 55455. Email:\,\texttt{njw@ima.umn.edu}.}, \, and\,
Shuzhong Zhang\footnote{Department of Systems Engineering and
Engineering Management, The Chinese University of Hong Kong,
Shatin, Hong Kong. Email:\,\texttt{zhang@se.cuhk.edu.hk}. Research
supported by Hong Kong RGC Earmarked Grants CUHK418505 and
CUHK418406.} } }

\maketitle

\begin{abstract}

In this paper we study the relationship between the optimal value
of a homogeneous quadratic optimization problem and that of its
Semidefinite Programming (SDP) relaxation. We consider two
quadratic optimization models: (1) $\min \{ x^* C x \mid x^* A_k x
\ge 1,\, x\in\mathbb{F}^n,\, k=0,1,...,m\}$; and (2) $\max \{ x^*
C x \mid x^* A_k x \le 1,\, x\in\mathbb{F}^n, \, k=0,1,...,m\}$.
If \emph{one} of $A_k$'s is indefinite while others and $C$ are
positive semidefinite, we prove that the ratio between the optimal
value of (1) and its SDP relaxation is upper bounded by $O(m^2)$
when $\mathbb{F}$ is the real line $\mathbb{R}$, and by $O(m)$
when $\mathbb{F}$ is the complex plane $\mathbb{C}$. This result
is an extension of the recent work of Luo {\em et
al.}~\cite{LSTZ}. For (2), we show that the same ratio is bounded
from below by $O(1/\log m)$ for both the real and complex case,
whenever all but one of $A_k$'s are positive semidefinite while
$C$ can be indefinite. This result improves the so-called
approximate S-Lemma of Ben-Tal {\em et al.}~\cite{BNR02}. We also
consider (2) with multiple indefinite quadratic constraints and
derive a general bound in terms of the problem data and the SDP
solution. Throughout the paper, we present examples showing that
all of our results are essentially tight.

%We consider two kinds of NP-hard quadratic optimization problems
%and study the performance of SDP relaxations to solve them.
%The first problem to minimize a convex quadratic form
%subject to $m$ concave homogeneous quadratic constraints plus
%one indefinite homogeneous quadratic constraint.
%We obtain the approximation bound
%$\mathcal{O}(m^2)$ for the real case and
%$\mathcal{O}(m)$ for the complex case,
%which is similar to the result in \cite{LSTZ}.
%The second problem to maximize a general (can be indefinite) quadratic form
%subject to $m$ convex homogeneous quadratic constraints plus
%one indefinite homogeneous quadratic constraint.
%We get the approximation bound
%$\mathcal{O}(1/\log m^2)$
%which is similar to the result in \cite{NRT}.
\end{abstract}

\bigskip


\medskip


\newpage

\section{Introduction}

%Before discussion, let us first start with the notations. We shall
%denote $\F$ to be either the field of real numbers $\re$ or
%complex numbers $\cpx$. For a matrix $A\in \F^{m\times n}$, let
%$A^*$ denote the conjugate transpose of $A$ (viz.\ $A^*=A^T$ if
%$\F=\re$, and $A^*=A^H$ if $\F=\cpx$). Moreover, let us denote
%$\SF^n$ to be the linear subspace of all (Hermitian) symmetric
%matrices in $\F^{n\times n}$, and $\SF^n_+$ to be the cone of
%positive semidefinite matrices in $\SF^n$, and $\SF^n_{++}$ to be
%the cone of positive definite matrices. Now, the problem that

We consider in this paper homogeneous quadratic optimization problems
in either the minimization form
\begin{equation}\label{minqp}
\begin{array}{rl}
\min &  x^* C x   \\
\st  &  x^* A_k x \geq 1, \, k=0,1,...,m \\
     &  x \in \F^n,
\end{array}
\end{equation}
or the maximization form
\begin{equation}\label{maxqp}
\begin{array}{rl}
\max &  x^* C x   \\
\st  &  x^* A_k x \leq 1, \, k=0,1,...,m \\
     &  x \in \F^n,
\end{array}
\end{equation}
where matrices $A_k$ and $C$ are $n\times n$, $\F$ can be the
field of real numbers $\re$ or the field of complex numbers
$\cpx$, and the superscript $^*$ represents Hermitian transpose
(or regular transpose in case of real numbers). Both of above
quadratic optimization problems are NP-hard \cite{LSTZ,BNR02},
even when all the data matrices, $C$ and $A_k$, $k=1,...,m$, are
positive semidefinite. Homogeneous quadratic optimization problems
(\ref{minqp})--(\ref{maxqp}) arise naturally in telecommunications
and robust control applications; see \cite{LSTZ,BNR02} and the
references therein. A popular approach to approximately solving
the NP-hard quadratic programs \reff{minqp}--\reff{maxqp} is to
use the so-called {\em Semidefinite Programming}\/ (SDP)
relaxations as follows:
\[
\begin{array}{rl}
\min &  \Tr (C X)   \\
\st  &  \Tr (A_k X) \geq 1, \, k=0,1,...,m \\
     &  X \in \SF^n_+,
\end{array}
\]
and, respectively,
\[
\begin{array}{rl}
\max &  \Tr (C X)    \\
\st  &  \Tr (A_k X) \leq 1, \, k=0,1,...,m \\
     &  X \in \SF^n_+,
\end{array}
\]
where $\Tr(\cdot)$ represents the trace of a matrix, $\SF^n_+$ denotes the convex cone of
positive semidefinite matrices in the space of all (Hermitian) symmetric
matrices $\SF^n$. %and $\SF^n_{++}$ to be
%the cone of positive definite matrices.
The above two SDPs are convex and can be solved efficiently via interior
point methods. After the SDP relaxation problems are solved, we
can apply a {\em randomization procedure} to the corresponding optimal SDP solutions to
extract rank-one feasible solutions for \reff{minqp} and \reff{maxqp}
respectively. Theoretically, even though the randomized solutions obtained in this
manner are not globally optimal for either
\reff{minqp} or \reff{maxqp}, they can be shown to be
high quality approximate solutions; see, e.g.~\cite{BNR02,LSTZ,NRT}.
%However, the quality of the randomized solution obtained this way
%depends critically on the tightness of the SDP relaxation.
Specifically,
Nemirovski {\em et al.}~\cite{NRT} proved that for the maximization
problem \reff{maxqp}, if all $A_k$'s are positive semidefinite, then the ratio
between the optimal value of the SDP relaxation problem and that
of the original quadratic problem is bounded above by $O(\log m)$.
More generally, Ben-Tal {\em et al.}~\cite{BNR02} established
a so-called {\em approximate S-Lemma} which
shows that the approximation ratio for the SDP
relaxation is at most $O(\log ( n^2 m ))$ when all but one of the
matrices $A_k$, $k=0,1,...,m$ are positive semidefinite.

%We note here that the
%approximate S-Lemma of Ben-Tal {\em et al.}~\cite{BNR02} has wide
%applications in robust optimization.
In a parallel development, Luo {\em et al.}~\cite{LSTZ}
considered the homogeneous quadratic optimization in minimization form
\reff{minqp}. It turns out
that the SDP approximation ratio for the
minimization version of the problem takes a quite different
form as compared with its maximization counterpart.
When all the matrices $A_k$ and $C$ are positive
semidefinite, Luo {\em et al.}~\cite{LSTZ} showed that the ratio
between the original optimal value and the SDP relaxation optimal
value is bounded above by $O(m^2)$ when $\F=\re$ and by $O(m)$
when $\F=\cpx$. All these bounds are shown to be tight in the
worst case, although the simulation studies in~\cite{LSTZ} showed
that the ratios are typically close to 1. In other words, the average
performance can be much better than the stated worst-case bounds for
randomly generated instances. %However, the worst-case performance
%analysis is still an extremely powerful indication of whether the
%method will work well in practice or not.
Recently, So {\em et al.}~\cite{SYZ06} developed methods for
finding approximate low rank solutions for linear matrix
inequalities. Their results unify the approximation bounds of
Nemirovski {\em et al.}~\cite{NRT} and Luo {\em et
al.}~\cite{LSTZ} as special cases (rank being 1), when all the
data matrices are positive semidefinite.

In this paper, we study the approximation ratio of the SDP
relaxation for homogeneous quadratic optimization problems
 \reff{minqp}--\reff{maxqp} when some of the constraint matrices $\{A_k\}$
 are {\em indefinite}.
Our results are as follows. In Section~\ref{HQP min form},  we show
that, for the problem in minimization form \reff{minqp},
the upper bounds for the approximation ratios of the SDP
relaxation as presented in~\cite{LSTZ} ($O(m^2)$ and $O(m)$
for $\F=\re$ and $\F=\cpx$ respectively)
hold true even when \emph{one} of the constraint matrices is indefinite.
If there are more than one indefinite quadratic
constraints, we show by an example that the approximation
ratio can be infinite. Therefore, our bounds are essentially best possible.
%Furthermore, we show in Section~\ref{GHQP} that if there are more
%than one indefinite quadratic constraints, then this ratio can
%even be infinite.
In Section~\ref{S-lemma}, we consider the problem in maximization form \reff{maxqp}.
We improve the approximate S-Lemma of Ben-Tal {\em et al.}~\cite{BNR02} by
reducing their upper bound on approximation ratio
from $O(\log ( n^2 m))$ to $O(\log m)$
when \emph{one} quadratic inequality is indefinite.  In the process of
establishing this new bound, we resolve
a conjecture by Ben-Tal {\em et al.}~\cite{BNR02}
on a possible universal lower bound for the probability that a homogeneous quadratic
form of binary i.i.d.\ Bernoulli random variables lies below its mean.
Finally, in Section~\ref{GHQP} we present a new and unifying upper bound on
the ratio of the optimal value of SDP relaxation over that of the original quadratic
maximization problem \reff{maxqp} {\em without}\/ any definiteness assumptions.
This new general bound involves the problem data and the SDP optimal
solution, which are computable in polynomial time. We also present
an example showing that this bound is essentially tight.

%However, if only one of $A_1,\cdots, A_m$ is indefinite, we can
%still prove the quality bound $ v_{qp} \leq \,\mathcal{O}(m^2) \,
%v_{sdp} $ for $\F=\re$, and $ v_{qp} \leq \,\mathcal{O}(m) \,
%v_{sdp} $ for $\F=\cpx$. The proof will be given in Section~2.

%\bigskip

%Another frequently used quadratic optimization is in the form
%\begin{align*}
%p_{qp}:=
%\max_{x\in \re^n} & \quad x^TAx \\
%s.t. & \quad  x^TA_k x \leq 1,\, k=1,...,m.
%\end{align*}
%where $A,A_1,\cdots,A_m$ are all symmetric matrices. Generally
%this problem is also known to be NP-hard. Its SDP relaxation is
%\begin{align*}
%p_{sdp}:=
%\max & \quad A\bullet X \\
%s.t. & \quad  A_k \bullet X \leq 1,\, k=1,...,m \\
%& \quad  X \succeq 0.
%\end{align*}
%When all $A_k(k=1,...,m)$ are positive semidefinite and
%$A_1+\cdots+A_m \succ 0$ (but $A$ can be indefinite), \cite{NRT}
%showed that
%\[
%p_{qp}\leq p_{sdp} \leq 2\log(2m\mu) \, p_{qp}
%\]
%where $\mu = \min\{m,\max_i Rank(A_i)\} \leq m$. Now we consider
%the case that there is one more constraint $x^TA_0x \leq 1$ where
%$A_0$ is symmetric indefinite. Then we can still prove similar
%quality bound, i.e.,
%\[
%p_{qp}\leq p_{sdp} \leq 2\log(174m\mu) \, p_{qp}.
%\]
%The proof will be given in Section~3. When there are two or more
%constraints like $x^TBx \leq 1$ with $B$ symmetric indefinite,
%there is no quality bound, as shown by one counterexample in
%Section~3.

\section{Estimating Asymmetry of a Random Variable About its Mean} \label{distribution below mean}

\setcounter{equation}{0}

To facilitate the technical analysis in subsequent sections, we
establish in this section a bound on the probability for a
general random variable to be above (or
symmetrically, below) its mean value, using only
the high order moment information of the random variable.
This problem is of importance on
its own in statistics and probability theory.
The following lemma is a generalization
of Theorem~2.1 in~\cite{KS}.

\begin{lemma} \label{mom4pb1}
Suppose that a random variable $\Phi$ satisfies $\ex \Phi=0$,
$\Var(\Phi)=1$ and $\ex |\Phi|^t \leq  \tau$ for some $t>2$ and
$\tau>0$. Then $\prob\{\Phi\geq 0\} > 0.25\tau^{-\frac{2}{t-2}}$
and $\prob\{\Phi\leq 0\}>0.25\tau^{-\frac{2}{t-2}}$.
\end{lemma}
\begin{proof}
Let $p_1=\prob\{\Phi\geq 0\}$ and $p_2=\prob\{\Phi\leq0\}$. Also let
$Y_1=\max(\Phi,0)$ and $Y_2=-\min(\Phi,0)$. Since $\ex\Phi=0$, we
know $\ex Y_1-\ex Y_2=0$. Let $s:=\ex Y_1= \ex Y_2$. By H\"older's
inequality it follows that $(\ex Y_1^t)^{1/(t-1)} (\ex
Y_1)^{(t-2)/(t-1)}\ge \ex Y_1^2$ and $(\ex Y_2^t)^{1/(t-1)} (\ex
Y_2)^{(t-2)/(t-1)}\ge \ex Y_2^2$. Since $\ex Y_1^t+\ex Y_2^t =\ex
|\Phi|^t $, we have
\[
\tau  \ge \ex |\Phi|^t= \ex Y_1^t+ \ex Y_2^t \geq \frac{(\ex
Y_1^2)^{t-1}+(\ex Y_2^2)^{t-1}}{s^{t-2}}.
\]
Let $u=\ex Y_1^2 \in [0,1]$. Since $\ex Y_1^2+\ex Y_2^2 =\ex
\Phi^2 = \Var(\Phi)=1$, it follows that $s^{t-2}\ge
\frac{u^{t-1}+(1-u)^{t-1}}{\tau}$. On the other hand, by the
Cauchy-Schwartz inequality, we have
\[
s^2 = (\ex Y_1)^2 = (\ex (\sign(Y_1)Y_1))^2 \leq \ex
(\sign(Y_1)^2) \ex Y_1^2 \leq p_1 u
\]
which implies that
\begin{eqnarray*}
p_1 &\ge&
u^{-1}\left[\frac{u^{t-1}+(1-u)^{t-1}}{\tau}\right]^{\frac{2}{t-2}} \\
&=& \frac{\left( u^{t-1}+(1-u)^{t-1} \right)^{\frac{2}{t-2}}}{u}
\,
\tau^{-\frac{2}{t-2}} \\
&\ge & \left( u^{t-1}+(1-u)^{t-1} \right)^{\frac{2}{t-2}}
\tau^{-\frac{2}{t-2}} \\
&\ge& \left( 2 \left( \frac{1}{2}\right)^{t-1}
\right)^{\frac{2}{t-2}} \tau^{-\frac{2}{t-2}} \\
&=& 0.25 \tau^{-\frac{2}{t-2}},
\end{eqnarray*}
where the third inequality follows from the convexity of the
function  $u^{t-1}$ when $t>2$. Obviously, the equality can not hold
throughout. Therefore, $p_1 > 0.25 \tau^{-\frac{2}{t-2}}$. By
symmetry, we also have $p_2 > 0.25 \tau^{-\frac{2}{t-2}}$.
\end{proof}

In case $t=4$, Lemma~\ref{mom4pb1} asserts that $\prob\{\Phi\geq
0\} \ge \frac{1}{4\tau}$ and $\prob\{\Phi \leq 0\} \ge
\frac{1}{4\tau}$. However, in this particular case, this specific
bound can in fact be further sharpened.

\begin{lemma} \label{mom4pb}
Suppose that a random variable $\Phi$ satisfies $\ex \Phi=0$,
$\Var(\Phi)=1$ and $\ex \Phi^4 \leq  \tau$. Then $\prob\{\Phi\geq
0\} \ge \frac{2\sqrt{3}-3}{\tau}> \frac{9}{20\tau}$ and $\prob\{
\Phi\leq 0\} \ge \frac{2\sqrt{3}-3}{\tau}>\frac{9}{20\tau}$.
\end{lemma}
\begin{proof}
It follows from the proof in the Lemma~\ref{mom4pb1} that
\[
p_1\ge \frac{u^3+(1-u)^3}{\tau
u}=\left(\frac{1}{u}+3u-3\right)\frac{1}{\tau}\ge
\frac{2\sqrt{3}-3}{\tau}
> \frac{9}{20\tau}.
\]
By symmetry, $p_2> \frac{9}{20\tau}$ holds as well.
\end{proof}

\section{Homogenous Quadratic Minimization and SDP Relaxation} \label{HQP min form}

\setcounter{equation}{0}

Consider the homogeneous quadratic optimization
%\begin{align*}
\begin{equation} \label{hqp1}
\begin{array}{rl}
v^{\min}_{qp}:=
\min &  x^* C x   \\
\st &  x^* A_k x \geq 1, \, k=0,1,...,m \\
& x \in \F^n,
%\end{align*}
\end{array}
\end{equation}
where $C,A_1,A_2,...,A_m\in \SF^{n}$ are symmetric matrices. This problem is generally
NP-hard \cite{LSTZ}. A natural semidefinite programming (SDP)
relaxation to the above quadratic optimization problem is
%\begin{align*}
\begin{equation}\label{sdp1}
\begin{array}{rl}
v^{\min}_{sdp}:=
\min &  \Tr(CZ)  \\
\st  & \Tr(A_kZ) \geq 1, \, k=0,1,...,m \\
     & Z \in \SF^n_+.
%\end{align*}
\end{array}
\end{equation}
Obviously, the SDP relaxation provides a lower bound, i.e.,
$v^{\min}_{sdp}\leq v^{\min}_{qp}$. In the case $C=I_n$, and
$A_0,A_1,...,A_m$ are all positive semidefinite, Luo {\em et
al.}~\cite{LSTZ} proved that $v^{\min}_{qp} / v^{\min}_{sdp} \le
\frac{27(m+1)^2}{\pi}$ for $\F=\re$, and $v^{\min}_{qp} /
v^{\min}_{sdp} \le 8(m+1)$ for $\F=\cpx$. Moreover, when there are two
or more of $A_0,A_1,...,A_m$ are indefinite, there is in general no
data-independent upper bound on $v^{\min}_{qp} /
v^{\min}_{sdp}$, as shown by the following example~\cite{LSTZ}:
\begin{align*}
\min & \quad x_1^2+x_2^2 \\
\st & \quad x_1^2 \geq 1 \\
& \quad x_1^2+Mx_1x_2 \geq 1 \\
& \quad x_1^2-Mx_1x_2 \geq 1
\end{align*}
where $M>0$ is a constant. In the above example, $v^{\min}_{sdp} =
1$, and the last two constraints imply $x_1^2\ge M|x_1| |x_2|+1$
which, together with the first constraint $x_2^2\ge 1$, yield
$x_1^2\ge M|x_1|+1$ or, equivalently,
$|x_1|\ge(M+\sqrt{M^2+4})/2$. Therefore, $v^{\min}_{qp} \ge
1+\frac{1}{4}(M+\sqrt{M^2+4})^2$. That is, $v^{\min}_{qp} /
v^{\min}_{sdp} \ge 1+\frac{1}{4}(M+\sqrt{M^2+4})^2$, which can be
arbitrarily large, depending on the problem data $M>0$.

%\section{Approximation bounds for quadratic optimization
%on the intersection of exteriors of ellipsoids}
%\setcounter{equation}{0}

In this section, we consider the homogeneous quadratic
optimization \reff{hqp1} under the assumption that
$C,A_1,A_2,...,A_m\in \SF^{n}_+$ are positive semidefinite
while $A_0\in \SF^{n}$ can be indefinite.
%\begin{align}
%v_{qp}:=
%\min_{x\in \F^n} & \quad x^TAx \label{hqp1} \\
%s.t. &  \quad x^TA_ix \geq 1, \, i=1,\cdots,m \label{hqp2} \\
%& \quad x^TA_0x\geq 1 \label{hqp3}
%\end{align}
%\begin{equation} \label{hqp1}
%\begin{array}{rl}
%v^{\min}_{qp}:= \min & x^* C x  \\
%\st                  & x^* A_k x \geq 1, \, k=1,...,m  \\
%& x^* A_0 x\geq 1 \\
%& x\in \F^n,
%\end{array}
%\end{equation}
%where $C,A_1,A_2,...,A_m\in \SF^{n}_+$ and $A_0\in \SF^{n}$, with
%the SDP relaxation
%%\begin{align}
%\begin{equation} \label{sdp1}
%\begin{array}{rl}
%v^{\min}_{sdp}:=
%\min & \Tr(CZ) \\ %\label{sdp1} \\
%\st  & \Tr(A_kZ) \geq 1, \, k=1,...,m \\ % \label{sdp2} \\
%&  \Tr(A_0Z) \geq 1 \\ %\label{sdp3} \\
%& Z \in \SF^n_+.
%\end{array}
%\end{equation}
%%\end{align}
Throughout this section, we assume that \reff{hqp1} is feasible, and that there is
$\mu_k\ge 0$, $k=0,1,...,m$, such that $ \sum_{k=0}^m
\mu_k A_k \prec 0$. This assumption guarantees that the SDP
relaxation is primal feasible while its dual problem satisfies the
Slater condition. Hence the strong duality holds and the primal problem
\reff{sdp1} has an optimal
solution that attains its infimum.

Our analysis shall treat the cases $\F=\re$ and $\F=\cpx$
separately, leading to different bounds and flavors. For clarity,
the analysis will be presented in the next two subsections.

\subsection{The real case} \label{real case}

Let us start with a useful lemma regarding a lower bound on worst
asymmetric mass distributions for a $\chi^2$-distribution around
its mean vector. In fact this result is interesting on its own
right.

\begin{lemma} \label{realpb}
Let $\tau_i $ be any real numbers, $i=1,...,n$, and let $\eta \sim
\,N(0,I_n)$ be an $n$-dimensional normal distribution with zero
mean and covariance matrix $I_n$. Then we have
\[
\prob \left\{ \sum_{i=1}^n \tau_i (\eta_i^2-1) \geq 0 \right\}
> \frac{3}{100}, \,\quad
\prob\left\{\sum_{i=1}^n \tau_i (\eta_i^2-1) \leq 0 \right\}
> \frac{3}{100}.
\]
\end{lemma}
\begin{proof}
Note that $\ex(\eta_i^2-1)^2 = \ex(\eta_i^4-2\eta_i^2+1) =
3-2+1=2$. Let $\Psi = \sum_{i=1}^n \tau_i (\eta_i^2-1)$, and $\Phi
=\frac{\Psi}{\sqrt{2\sum_{i=1}^n \tau_i^2}}$. Then $\ex \Phi=0$
and $\Var(\Phi)=1$. Since $ \ex (\eta_i^2-1)^2 =2$, and $\ex
(\eta_i^2-1)^4 = 60$, direct calculation shows
\[
\ex \Psi^4 =  48 \sum_{i=1}^n \tau_i^4 + 12 \left(\sum_{i=1}^n
\tau_i^2\right)^2 \leq 60 \left(\sum_{i=1}^n \tau_i^2\right)^2.
\]
Therefore, we have
\[
\ex \Phi^4 = \frac{\ex \Psi^4}{4(\sum_{i=1}^n \tau_i^2)^2} \leq
15.
\]
It follows from Lemma~\ref{mom4pb} that $\prob \{ \Phi \geq 0 \} >
\frac{3}{100}$. Similarly, we have $\prob\{\Phi \leq 0 \} >
\frac{3}{100}$ by symmetry.
\end{proof}

%\begin{remark}
Using H\"older's inequality, we also have $\ex |\Psi|^3 \leq
60^{\frac{3}{4}} (\sum_{i=1}^n \tau_i^2)^{\frac{3}{2}}$ and $\ex
|\Phi|^3 \leq 15^{\frac{3}{4}}$ which can be used to
lower $\prob \{ \Phi
\geq 0 \}$ (c.f.\ Theorem~2.1 in \cite{KS}).
However, in this particular case, the bound so obtained is
slightly worse than the one that we derived in Lemma~\ref{realpb}.
%\end{remark}

\begin{lemma}\label{newlemma}
Let $A,Z$ be two real symmetric matrices with $Z \succeq 0$ and $
\Tr(A Z) \ge 0$. Let $\xi \in N(0,Z)$ be a normal random vector
with zero mean and covariance matrix $Z$. Then for any $0\le
\gamma \le 1 $ we have
\begin{align*}
\prob\{ \xi^T A \xi < \gamma \ex (\xi^T A \xi) \} <\, 1  -
\frac{3}{100} .
\end{align*}
\end{lemma}
\begin{proof}
Let $r=\rank(AZ)$, and $Q\in \re^{n\times n}$ be an orthogonal
matrix such that
\[
Q^T(Z^{\half} A Z^{\half})Q = \text{diag}(\lmd_1,\cdots, \lmd_r,
0, \cdots ,0) .
\]
Since $\Tr(A Z)\ge 0$ we have $\sum_{i=1}^r \lmd_i \ge 0$. Let
$\bar\xi \in N(0,I_n)$ and $\xi := Z^\half Q \bar\xi$. Then
$\xi$ follows a Gaussian distribution $ N(0,Z)$.
%Then
%\[
%\bar \xi = Q^T(Z^\half)^+\xi \in N\left(0, \left[\begin{array}{cc}
%I_r & 0_{r\times (n-r)} \\ 0_{(n-r)\times r} & 0_{(n-r)\times
%(n-r)}
%\end{array} \right] \right)
%\]
%is a normal random vector, where $(Z^\half)^+$ stands for the
%pseudo-inverse of $Z^\half$.
Moreover, we have $ \xi^T A \xi = \sum_{i=1}^r \lmd_i \bar\xi_i^2
$, where $\bar\xi_i$, $i=1,...,r$, are independent and follow the
normal distribution $N(0,1)$. Therefore, we have $\ex(\xi^T A
\xi)=\sum_{i=1}^r \lmd_i$ and
%\begin{align*}
\begin{eqnarray*}
\prob \{ \xi^T A \xi < \gamma \ex (\xi^T A \xi) \}
& = & \prob \left\{ \sum_{i=1}^r \lmd_i \bar\xi_i^2 < \gamma \sum_{i=1}^r \lmd_i\right\} \\
& = & \prob \left\{\sum_{i=1}^r \lmd_i (\bar\xi_i^2-1) < (\gamma-1) \sum_{i=1}^r \lmd_i\right\} \\
& \leq& \prob \left\{ \sum_{i=1}^r \lmd_i (\bar\xi_i^2-1) <
0\right\}  <  1 - \frac{3}{100},
\end{eqnarray*}
%\end{align*}
where the first inequality follows from $\gamma\in[0,1]$ and
$\sum_{i=1}^r \lmd_i\ge 0$, and the
last step is due to Lemma~\ref{realpb}.
\end{proof}

Now we are ready to establish the following quality bound for the SDP relaxation.
The argument follows closely those of \cite{LSTZ}.

\begin{thm} \label{real-case Prob anal}
Consider the real quadratic program $\reff{hqp1}$ and its SDP relaxation
$\reff{sdp1}$, where $\F=\re$. Then, there holds
\[
\frac{v^{\min}_{qp}}{v^{\min}_{sdp}} \leq \frac{10^6m^2}{\pi}.
\]
\end{thm}
\begin{proof}
Let $\hat Z$ be an optimal solution of the SDP relaxation
\reff{sdp1} with rank $r$ satisfying $\frac{(r+1)r}{2}\leq m$. The
existence of such matrix solution is well known;
cf.~Pataki~\cite{Pata}. Moreover, this low rank matrix can be
constructed in polynomial-time; cf.~\cite{HZ}. Clearly, $r <
\sqrt{2m}$. Since $\hat Z$ is
feasible, $\Tr(A_0\hat Z) \geq 1$. %Applying previous lemmas, we get
For any $0<\gamma\le 1$ and $\mu>0$ we have
%\begin{align*}
\begin{eqnarray*}
& & \prob \left\{ \min_{0\leq k\leq m} \xi^TA_k\xi \geq \gamma,\,  \xi^TC\xi \leq \mu \, \Tr(C\hat Z)\right\} \\
&=& \prob \left\{ \xi^TA_k\xi \geq \gamma \, \mbox{ for all }k=0,1,...,m,
\text{ and }
\,  \xi^T C \xi  \leq \mu \, \Tr(C\hat Z)\right\}  \\
&\geq & \prob  \left\{ \xi^TA_k\xi \geq \gamma\, \Tr(A_k\hat Z) \,
\mbox{ for all }k=0,1,...,m, \text{ and }
\,  \xi^T C \xi  \leq \mu \, \Tr(C\hat Z)\right\} \\
&=& \prob \left\{\xi^T A_k \xi \geq \gamma\, \ex (\xi A_k \xi) \,
\mbox{ for all }k=0,1,...,m, \text{ and }
\,  \xi^T C\xi  \leq \mu \,  \ex (\xi^TC\xi) \right\} \\
& \geq & 1-\sum_{k=0}^m \prob  \left\{ \xi^TA_k\xi <\gamma\, \ex
(\xi A_k \xi) \right\} - \prob \left\{ \xi^T C\xi > \mu \, \ex
(\xi^TC\xi) \right\}.
\end{eqnarray*}
%\end{align*}
Since $A_k \succeq 0$ for $k=1,...,m$, it follows from Lemma~3.1
of~\cite{LSTZ} that
\[
\prob \left\{ \xi^T A_k \xi < \gamma \ex (\xi^TA_k\xi) \right\}
\leq
 \max\left\{\sqrt{\gamma}, \frac{2(r-1)\gamma}{\pi-2} \right\}.
\]
Although $A_0$ is indefinite,  we can use Lemma~\ref{newlemma} to
obtain
\[
\prob \left\{ \xi^T A_0 \xi < \gamma \ex (\xi^TA_0 \xi) \right\} <
1-\frac{3}{100}.
\]
Also, since $C\succeq 0$, we can apply Markov inequality to obtain
$$ \prob
\left\{\xi^T C\xi > \mu \, \ex (\xi^TC\xi) \right\} \leq
\frac{1}{\mu}. $$
Combining the above estimates yields
%\begin{align*}
%\begin{eqnarray*}
\[
 \prob \left\{ \min_{0\leq k\leq m} \xi^TA_k\xi \geq \gamma ,\,
 \xi^T C \xi \leq \mu \, \Tr(C\hat Z)\right\}
 >  \frac{3}{100} -m \max\left\{\sqrt{\gamma},
\frac{2(r-1)\gamma}{\pi-2} \right\} -\frac{1}{\mu}.
\]
%\end{eqnarray*}
%\end{align*}

Let $ \mu = 100$ and $\gamma =\frac{\pi}{10^4m^2} $. %, and so $\gamma<1$, and
Since $r <\sqrt{2m}$, we have $\sqrt{\gamma} \geq
\frac{2(r-1)\gamma}{\pi-2}$. %, \forall \, m=1,2, \cdots
%\]
For these values of $\mu$ and $\gamma$, we have
\[
\frac{3}{100} -m \max\left\{\sqrt{\gamma},
\frac{2(r-1)\gamma}{\pi-2} \right\} -\frac{1}{\mu} = \frac{3}{100}
-m \frac{\sqrt{\pi}}{100m} - \frac{1}{100} > \frac{1}{500} .
\]
Therefore, there exists a vector $\xi \in \re^n$ such that
\[
\xi^TA_k\xi \geq \gamma,\,\, \ k=0,1,...,m,\, \ \mbox{ and } \ \xi^TC\xi \leq \mu
\Tr(C\hat Z).
\]
Now let $x =  \frac{1}{\sqrt{\gamma}} \xi$. Then, $x^TA_k x\geq
1$, $k=0,1,...,m$, and
\[
v^{\min}_{qp} \leq x^TCx = \frac{1}{\gamma} \xi^TC\xi \leq
\frac{\mu}{\gamma} \Tr(C\hat Z) = \frac{10^6m^2}{\pi}
\,v^{\min}_{sdp},
\]
which establishes the desired bound
\end{proof}

\subsection{The complex case} \label{complex case}

Recall that the density function of a complex-valued normal
distribution\footnote{For a discussion on the complex normal
distribution and the related references, see Zhang and
Huang~\cite{ZH06}.} $\eta \sim N_c(0,1)$ is
\[
\frac{1}{\pi}e^{-|u|^2},\, \forall u \in \cpx.
\]
In polar coordinates, the density function becomes
\[
\frac{\rho}{\pi}e^{-\rho^2},\, \forall\,\rho\in [0,+\infty),\, \theta \in [0,2\pi).
\]
The argument $\theta$ is uniformly distributed in $[0,2\pi)$, and
the modulus $\rho$ has the distribution
\[
f(\rho) =
\bca
2\rho e^{-\rho^2}, & \text{if } \rho \geq 0 ; \\
0, & \text{if } \rho <0.
\eca
\]
Thus squared modulus $|\eta|^2$ has the exponential distribution
\[
\prob \{ |\eta|^2 \leq \alpha \} \leq 1  - e^{-\alpha}.
\]

\begin{lemma} \label{cpxpb}
For any real numbers $\tau_i $, and i.i.d.\ exponential random variables
$\eta_i$ with unit variance, $i=1,...,n$, there holds
\[
\prob \left\{\sum_{i=1}^n \tau_i (\eta_i-1) \geq 0 \right\}
> \frac{1}{20},\qquad
\prob\left\{ \sum_{i=1}^n \tau_i (\eta_i-1) \leq 0 \right\}
> \frac{1}{20}.
\]
\end{lemma}
\begin{proof}
Note that $ \ex(\eta_i-1)^2 = 1$. Let $\Psi = \sum_{i=1}^n \tau_i
(\eta_i-1)$ and $\Phi =\frac{\Psi}{\sqrt{\sum_{i=1}^n \tau_i^2}}$.
Clearly, $\ex \Phi=0$ and $\Var(\Phi)=1$. Since $\ex (\eta_i-1)^4
= 9$, direct calculation shows
\[
\ex \Psi^4 =  6 \sum_{i=1}^n \tau_i^4 + 3 \left(\sum_{i=1}^n
\tau_i^2\right)^2 \leq 9 \left(\sum_{i=1}^n \tau_i^2\right)^2.
\]
This further implies
\[
\ex \Phi^4 = \frac{\ex \Psi^4}{(\sum_{i=1}^n \tau_i^2)^2} \leq 9.
\]
Using Lemma~\ref{mom4pb} we have $\prob\{ \Phi \geq 0 \} >
\frac{1}{20}$. Similarly, $\prob\{ \Phi \leq 0 \}
> \frac{1}{20}$.
\end{proof}

\bigskip

%\begin{conj}
%For nonnegative $\tau_i $ such that $\sum_{i=1}^k \tau_i =1$
%and random variable $\eta_i$ with exponential distribution, it holds
%\[
%\text{Prob}\left(\sum_{i=1}^k \tau_i (\eta_i-1) \geq 0 \right)
%\geq e^{-1}.
%\]
%\end{conj}
Interestingly, it is possible to find a closed formula (see
e.g.~\cite{Cox} and~\cite{Amari}) for the above probability. In
particular, if all the $\tau_i$'s are distinctive, then
\[
\prob \left\{\sum_{i=1}^n \tau_i (\eta_i-1) \geq 0 \right\} =
\sum_{i=1}^n \frac{e^{-\frac{1}{\tau_i}}} {\prod_{j\ne i}
\left(1-\frac{\tau_j}{\tau_i}\right)}.
\]
Therefore, we have
\[
\frac{1}{20} < \sum_{i=1}^n \frac{e^{-\frac{1}{\tau_i}}}
{\prod_{j\ne i} \left(1-\frac{\tau_j}{\tau_i}\right)} <
\frac{19}{20}
\]
for any distinctive real values $\tau_i$, $i=1,...,n$.

In fact, we conjecture that the following tighter inequalities
\begin{equation} \label{ineq conj}
\frac{1}{e} < \sum_{i=1}^n \frac{e^{-\frac{1}{\tau_i}}}
{\prod_{j\ne i} \left(1-\frac{\tau_j}{\tau_i}\right)} <
\frac{e-1}{e},
\end{equation}
hold for any real values $\tau_i$, $i=1,...,n$.
Inequality \reff{ineq conj} can be shown to hold for $n=2,3$. It also
admits a geometric interpretation. Specifically, let us consider the joint
exponential distribution on $\re_+^n$ with density $e^{-\sum_{i=1}^nx_i}$.
Then, the mean vector of this distribution, or equivalently,
the center of gravity of $\re_+^n$ is $x^c:=(1,1,...,1)^T$.
Given any real numbers $\tau_i$, $i=1,...,n$, the set
\[
{\cal H}= \left\{(\eta_1,\eta_2,...,\eta_n)^T\, \left|\,
\sum_{i=1}^n \tau_i (\eta_i-1) = 0\right\} \right.
\]
represents a hyperplane passing through $x^c$. If we let ${\cal H}_+$ denote the half space
in $\re^n$ created by the positive side of ${\cal H}$, then
inequality \reff{ineq conj} can be
interpreted as follows:
$$\mbox{Prob}(\re_+^n\cap {\cal H}_+)\ge e^{-1},\quad
\mbox{for any hyperplane ${\cal H}$ passing through $x^c$.}$$
Interestingly, the well-known theorem of
Gr\"{u}nbaum~\cite{Grun} can also be viewed from
this perspective: for any bounded convex body $\mathcal{C}\subset
\re^n$, if we assign the uniform distribution to $\mathcal{C}$, then
the mean vector of this distribution is given by the center of gravity
\[
x^c=\frac{1}{\mbox{Volume}(\mathcal{C})}\int_{\mathcal{C}}dx;
\]
as a result, if we consider any hyperplane $\cal H$ passing through
$x^c$ and let ${\cal H}_+$ denote the positive side of the hyperplane, then
Gr\"{u}nbaum inequality
\[
\mbox{Volume}\left(\mathcal{C}\cap \mathcal{H}_+\right) \geq
\,e^{-1}\,\, \mbox{Volume}\left(\mathcal{C}\right)
\]
can be written as
$$\mbox{Prob}(\mathcal{C}\cap\mathcal{H}_+)\ge e^{-1},\quad
\mbox{for any hyperplane $\cal H$ passing through $x^c$.}$$
Thus, inequality \reff{ineq conj} can be viewed as an extension of
Gr\"{u}nbaum's theorem to the exponential distribution over the unbounded
convex set $\mathcal{C}=\re_+^n$.

\bigskip

\begin{lemma} \label{complex case lemma}
Let $A,Z$ be two Hermitian matrices satisfying $ Z \succeq 0$ and
$\Tr(A Z) \ge 0$. Let $\xi \sim N_c(0,Z)$ be a complex normal
random vector. Then, for any $0\le \gamma \le 1$, we have
%\begin{align*}
\[
\prob\{\xi^* A \xi < \gamma \ex (\xi^* A \xi) \} < 1  -
\frac{1}{20}.
\]
%\end{align*}
\end{lemma}
\begin{proof}
Let $Q\in \cpx^{n\times n}$ be an unitary matrix such that
\[
Q^*(Z^\half A Z^\half)Q = \text{diag}(\lmd_1,\cdots,
\lmd_r,0,\cdots,0)
\]
where $r=\rank(A Z)$. Since $\Tr(A Z)\ge 0$, it follows that
$ \sum_{i=1}^r \lmd_i
\ge 0$. Let $\hat \xi\in\cpx^n$ be a random Gaussian vector drawn from
the complex normal distribution $N_c(0,I_n)$. Then the random vector
$\xi = Z^\half Q \hat \xi$ follows the Gaussian distribution $N_c(0,Z)$.
As a result, there holds
%\begin{align*}
\begin{eqnarray*}
\prob \left\{ \xi^* A \xi < \gamma \ex (\xi^* A \xi)\right)
&=& \prob \left\{\sum_{i=1}^r \lmd_i |\hat \xi_i|^2
< \gamma  \sum_{i=1}^n \lmd_i   \right\} \\
&=& \prob \left\{\sum_{i=1}^n \lmd_i  (|\hat \xi_i|^2-1) <
(\gamma-1)  \sum_{i=1}^n \lmd_i   \right\} \\
& \leq& \prob \left\{ \sum_{i=1}^n \lmd_i  (|\hat \xi_i|^2-1) < 0
\right \},
\end{eqnarray*}
%\end{align*}
where the last step follows from $\gamma\in[0,1]$ and $ \sum_{i=1}^r \lmd_i
\ge 0$.
Since $|\xi_i|^2$ is exponentially distributed, by
Lemma~\ref{cpxpb}, we have
\[
\prob \left\{\sum_{i=1}^n \lmd_i (|\hat \xi_i|^2-1) \geq 0
\right\}
> \frac{1}{20}
\]
which proves the lemma.
\end{proof}

\begin{thm}\label{thm2}
Consider $\reff{hqp1}$ and $\reff{sdp1}$, where $\F=\cpx$. Then
\[
\frac{v^{\min}_{qp}}{v^{\min}_{sdp}} \leq 2400 m .
\]
\end{thm}
\begin{proof}
It is known that in this case, if $v^{\min}_{sdp}$ is finite and
$m\le 3$, then $v^{\min}_{qp}/v^{\min}_{sdp}=1$
(cf.~e.g.~\cite{HZ} and~\cite{YZ03}). Below we shall only consider
the case where $m\ge 4$. Let $\hat Z$ be a low rank optimal
solution of the SDP relaxation \reff{sdp1}, such that
$r=\rank(\hat Z) \le \sqrt{m}$ (see~\cite{HZ}, \S5).
%and $Z^*$ can
%be chosen to have rank $r=1$ when $m\leq 3$ and $r\leq \sqrt{m}$
%when $m\geq 4$; see \cite[\S5]{HZ}. So when $m\leq 3$, the SDP
%relaxation is exact, and thus $v_{qp} = v_{sdp}$.
%Now consider the case $m\geq 4$.
The feasibility of $\hat Z$ implies that $\Tr(A_0 \hat Z) \ge
1$. Similar to Theorem~\ref{real-case Prob anal}, we can use the union
bound to obtain the following inequality %Now we estimate the following probability
%\begin{align*}
\begin{eqnarray*}
&&\prob \left\{ \min_{0\leq k\leq m} \xi^* A_k\xi \geq \gamma , \, \xi^* C\xi \leq \mu \, \Tr(C\hat Z)\right\} \\
%& = & \prob \left\{ \xi^* A_k\xi \geq \gamma \, \mbox{ for all }
%k=0,1,...,m,\ \text{and}\
%\, \xi^* C \xi \leq \mu \, \Tr(C \hat Z) \right\} \\
%& \geq & \prob\left\{ \xi^* A_k\xi \geq \gamma\, \Tr(A_k \hat Z)
%\, \mbox{ for all }k=0,1,...,m,\ \text{and}\
%\, \xi^* C\xi \leq \mu \, \Tr(C\hat Z)\right\} \\
%& = & \prob \left\{ \xi^* A_k \xi \geq \gamma\, \ex (\xi^* A_k
%\xi) \, \mbox{ for all }k=0,1,...,m, \ \text{and}\
%\, \xi^* C \xi \leq \mu \,  \ex (\xi^* C \xi) \right\} \\
& \geq & 1-\sum_{k=0}^m \prob \left\{ \xi^* A_k \xi <\gamma\, \ex
(\xi^* A_k \xi) \right\} - \prob \left\{ \xi^*C\xi > \mu \, \ex
(\xi^*C\xi) \right\}.
\end{eqnarray*}
%\end{align*}
Since $A_k\succeq 0$, $k=1,...,m$, it follows from Lemma~3.4 in~\cite{LSTZ} that
\[
\prob \left\{ \xi^* A_k \xi < \gamma \ex (\xi^* A_k\xi) \right\}
\leq \max\left\{\frac{4}{3}\gamma, 16(r-1)^2\gamma^2\right\}.
\]
Although $A_0$ is indefinite, Lemma~\ref{complex case lemma},
 asserts that
\[
\prob \left\{ \xi^* A_0 \xi < \gamma \ex (\xi^* A_0 \xi) \right\}
< 1-\frac{1}{20}.
\]
Therefore, combining these estimates and using Markov inequality, we have
%\begin{align*}
\[
\prob \left\{ \min_{0\leq k\leq m} \xi^* A_k\xi \geq \gamma ,
\xi^* C\xi \leq \mu, \, \Tr(C\hat Z)\right\}
> \frac{1}{20}
-m\max\left\{\frac{4}{3}\gamma, 16(r-1)^2\gamma^2\right\}
-\frac{1}{\mu}.
\]
%\end{align*}
Now choose $ \mu = 60$ and $\gamma =\frac{1}{40m}$. In this case,
%makes $\gamma<1$ and
$\frac{4}{3}\gamma \geq 16(r-1)^2\gamma^2$. We also have a strict
lower bound of the above probability
\[
\prob \left\{ \min_{0\leq k\leq m} \xi^* A_k\xi \geq \gamma , \,
\xi^* C\xi \leq \mu \, \Tr(C\hat Z)\right\} > 0.
\]
This implies that there exists $\xi\in \cpx^n$ such that
\[
\xi^* A_k \xi \geq \gamma, \, k=0,1,...,m;\quad \xi^* C \xi \leq \mu
\Tr(C \hat Z).
\]
Now let $x :=  \frac{1}{\sqrt{\gamma}} \xi$. Then $x^* A_k x\geq
1$, $k=0,1,...,m$, and so
\[
v^{\min}_{qp} \leq x^* C x \leq \frac{\xi^* C \xi}{\gamma} \leq
\frac{\mu\Tr (C \hat Z)}{\gamma} =  2400 m\cdot v^{\min}_{sdp}.
\]
The theorem is proven.
\end{proof}

Notice that there are examples (see \cite{LSTZ}) which show that the worst-case
ratios of $v^{\min}_{qp} /v^{\min}_{sdp} $ are indeed $O(m^2)$ and $O(m)$
in the real and complex case respectively, even in the absence of indefinite
constraint $x^*A_0x\ge 1$. Thus, the bounds of Theorems~\ref{real-case Prob anal}
and \ref{thm2}
are essentially tight.

Finally, we may also
wonder what happens if there are more than one indefinite quadratic constraint.
The following example shows that in this case the SDP relaxation
does not admit {any} \emph{finite} quality
bound.

\begin{example} \label{QP min 4 constraints}
\[
\begin{array}{rl}
\min & x_4^2 \\
\st  & x_1 x_2 + x_3^2 + x_4^2 \ge 1 \\
     & -x_1 x_2 + x_3^2 + x_4^2 \ge 1 \\
     & \frac{1}{2} x_1^2 - x_3^2 \ge 1 \\
     & \frac{1}{2} x_2^2 - x_3^2 \ge 1 \\
     & x_1,x_2,x_3,x_4 \in \re.
\end{array}
\]
The first two constraints are equivalent to $|x_1x_2| \le
x_3^2+x_4^2-1$. %, and so
%\begin{equation} \label{lower side}
%x_1^2 x_2^2 \le (x_3^2+x_4^2-1)^2.
%\end{equation}
At the same time, the last two constraints imply
\(%\begin{equation} \label{upper side}
|x_1 x_2| \ge 2 (x_3^2+1).
\) %\end{equation}
Combining these two inequalities yields %\reff{lower side} with \reff{upper side} yields
\[
x_3^2+x_4^2-1\ge 2 (x_3^2+1),
\]
which further implies $x^2_4\ge 3$.
%\[
%2(x_3^4+x_4^4)+1 \ge (x_3^2+x_4^2)^2 - 2 (x_3^2+x_4^2) + 1 =
%(x_3^2+x_4^2-1)^2 \ge x_1^2 x_2^2 \ge 4 x_3^4 + 8 x_3^2 + 4.
%\]
%Thus, $2x_4^4 \ge 2 x_3^4 + 8 x_3^2 +3\ge 3$, and so $x_4^2 >1$.
Therefore, we must have $v^{\min}_{qp} \ge 3$ in this case. However,
\[
\left[ \begin{array}{cccc} 4 & 0 & 0 & 0 \\ 0 & 4 & 0 & 0 \\ 0 & 0
& 1 & 0 \\ 0 & 0 & 0 & 0 \end{array} \right]
\]
is feasible for the corresponding SDP relaxation problem and
attains an objective value of $0$. Thus, it must be
optimal and thus $v^{\min}_{sdp}=0$.
Hence, $v^{\min}_{qp}/v^{\min}_{sdp} = \infty$ in this case.
\end{example}

\section{Quadratic Maximization and the Approximate $S$-Lemma} \label{S-lemma}

%\section{Approximation bound on maximizing quadratic forms over
%intersection of ellipsoids and a hyperboloid }

\setcounter{equation}{0}

In this section, we consider the nonconvex homogeneous quadratic optimization in the
maximization form
%\begin{align}
\begin{equation} \label{Qp-max}
\begin{array}{rl}
v^{\max}_{qp}:=
\max & x^* C x \\ %\label{eqp1} \\
\st  & x^* A_k x \leq 1,\, k=0,1,...,m \\ % \label{eqp2}\\
%& x^* A_0 x \le 1 \\ % \label{eqp3} \\
& x \in \F^n,
\end{array}
\end{equation}
%\end{align}
where $A_k \in \SF^n_+$, $k=1,...,m$, are positive semidefinite,
while $C, A_0 \in \SF^n$ may be indefinite. For convenience, from now on we shall
focus on the case $\F=\re^n$. Unlike the case of minimization form, this choice does
not significantly affect the quality of SDP approximation ratios,
since in the complex case the bounds are of the same order of magnitude.
%The feasible region \reff{eqp2}-\reff{eqp3} is the intersection of $m$
%ellipsoids $\mathcal{E}_i=\{x\in \re^n:\, x^TA_ix\leq 1\}$ and a
%hyperboloid $\mathcal{H}=\{x\in \re^n:\, x^TA_0x\leq 1\}$, which
%is generally nonconvex. The objective function $x^TAx$ is also
%nonconvex. We assume matrix $A_0$ is positive definite on the
%subspace $\cap_{i=1}^m \mathcal{N}(A_i)$, which is equivalent to
%the following set
%\[
%\Gamma := \left\{0\ne x\in\re^n:
%x^TA_0x\leq 0, A_1x=0,\cdots,A_m x =0
%\right\}
%\]
%is empty. Otherwise, if $\Gamma \ne \emptyset$,
%then the feasible set \reff{eqp2}-\reff{eqp3}
%is unbounded and the optimal value might be
%infinity for some particular matrix $A$.
We assume that there is $\mu_k\ge 0$, $k=0,1,...,m$, such that
\[
\sum_{k=0}^m \mu_k A_k \succ 0.
\]
Under this condition, the SDP relaxation satisfies the dual
Slater condition. Thus the primal-dual optimal solutions exist
and the primal-dual optimal objective
values are attainable. Let the SDP relaxation optimal value be
%\begin{align}
\begin{equation}\label{sdp2}
\begin{array}{rl}
v^{\max}_{sdp}:=
\max &  \Tr(C  X) \\ % \label{esdp1}\\
\st  &  \Tr(A_k X) \leq 1,\, k=0,1,...,m \\ % \label{esdp2}\\
%&   A_0 \bullet X \le 1, \\ % \label{esdp3}\\
&   X \succeq 0. \\ % \label{esdp4}
%\end{align}
\end{array}
\end{equation}
Obviously $v^{\max}_{qp}\leq v^{\max}_{sdp}$.
%In this section we shall prove that $p_{sdp} \leq C p_{qp}$ for some constant $C$.

%implies the set $\Gamma\ne \emptyset$ which is not possible by our assumption.

\begin{lemma} \label{qformprob}
Let $w_{ij}$ $(1\leq i<j\leq n)$ be any real numbers, and $\xi_i$
$(1\leq i\leq n)$ be random variables such that $\prob\left\{
\xi_i=-1\right\}=\prob\left\{\xi_i=1\right\}=0.5$. Then there holds
\[
\prob \left\{ \sum_{1\leq i<j\leq n } w_{ij} \xi_i \xi_j \le 0
\right\} > \frac{1}{87}.
\]
\end{lemma}
\begin{proof}
Let $\Psi =\sum_{1\leq i<j\leq n } w_{ij} \xi_i \xi_j$. Then $\ex
\Psi=0$, $\ex (\Psi^2)=\sum_{1\leq i<j\leq n } w_{ij}^2$ and
\[
\ex (\Psi^4) = \sum_{1\leq i<j\leq n } w_{ij}^4  + 6\sum_{(i,j)<
(k,\ell)} w_{ij}^2 w_{k\ell}^2 + W
\]
where $(i,j)< (k,\ell)$ means $i<k$ or $i=k$ and $j<\ell$, and
%\begin{align*}
\begin{eqnarray*}
W &=& 24 \sum_{1\leq i<j<k<\ell\leq n}
\left(w_{ij}w_{ik}w_{j\ell}w_{k\ell}+
w_{ij}w_{i\ell}w_{jk}w_{k\ell}+w_{ik}w_{i\ell}w_{jk}w_{j\ell}\right) \\
& \leq & 6 \sum_{1\leq i<j<k<\ell\leq n} \left(
(w_{ij}^2+w_{k\ell}^2)(w_{ik}^2+w_{j\ell}^2)+
(w_{ij}^2+w_{k\ell}^2)(w_{i\ell}^2+w_{jk}^2)+
(w_{ik}^2+w_{j\ell}^2)(w_{i\ell}^2+w_{jk}^2) \right) \\
& \leq & 36 \left(\sum_{1\leq i<j \leq n} w_{ij}^2\right)^2.
\end{eqnarray*}
%\end{align*}
Therefore we have $\ex (\Psi^4) \leq 39 (\sum_{1\leq i<j\leq n }
w_{ij}^2)^2$, since
\[
\sum_{1\leq i<j\leq n } w_{ij}^4  + 6\sum_{(i,j)< (k,\ell)}
w_{ij}^2 w_{k\ell}^2  \leq 3 \left(\sum_{1\leq i<j \leq n} w_{ij}^2\right)^2.
\]
Now let $\Phi =\frac{\Psi}{\sqrt{\sum_{1\leq i<j\leq n }
w_{ij}^2}}$. Then $\ex (\Phi)=0$, $\Var(\Phi)=1$ and $ \ex
(\Phi^4) \leq 39$. By Lemma~\ref{mom4pb}, we have
\[
\prob \left\{ \Phi \leq 0 \right\} >  \frac{1}{87}.
\]
The desired result follows.
\end{proof}

Lemma~\ref{qformprob} settles in the affirmative an open
question of Ben-Tal {\em et al.}~\cite[Conjecture A.5]{BNR02} who
conjectured that
\[
\prob \left\{ \sum_{1\leq i<j\leq n } w_{ij} \xi_i \xi_j \le 0
\right\} \ge \frac{1}{4},\qquad \forall \ w_{ij},
\]
except that we have a smaller  constant of $1/87$.
The above inequality was needed to establish the so called
{\em approximate $S$-Lemma} --- an extension of the well-known
$S$-Lemma, which is important in the context of robust
optimization and is closely related to our analysis in this
section.
In their work \cite{NRT}, Ben-Tal {\em et al.}~derived a weaker
lower bound of $1/8n^2$, which goes to zero as $n\to\infty$.

We can now use Lemma~\ref{qformprob} to analyze the performance
of SDP relaxation for \reff{sdp2}.
Let $\hat X=UU^T$ be one optimal solution of \reff{sdp2}, where $U\in \re^{n\times
r}$ and $r=\rank(\hat X)$. Suppose $Q \in \re^{n\times r}$ is the
orthogonal matrix such that $\hat C:=Q^TU^TCUQ$ is diagonal. Let
$\xi_k$, $k=1,...,r$, be i.i.d.\ random variables taking values
$-1$ or $1$ with equal probabilities, and let
\[
x(\xi) := \frac{1}{\displaystyle \sqrt{\max_{0\leq k\leq m}  \xi^T\hat A_k\xi}}
\,\, UQ\xi,
\]
where $\hat A_k = Q^TU^TA_k UQ$. Note that the above random vector
$x(\xi)$ is always well-defined, since  the assumption  $\sum_{k=0}^m \mu_k A_k\succ 0$
implies
\[
\max_{0\leq k\leq m} \xi^T\hat A_k\xi > 0 \mbox{  for any }
\xi\not=0.
\]
Let $\mu = \min\{m,\max_i \rank(A_i\hat X)\}$. We have the
following estimate of the SDP approximation ratio.
\begin{thm} \label{max qp bound}
There holds
\[
v^{\max}_{qp} \leq v^{\max}_{sdp} \leq 2\log(174\,m\mu)\,
v^{\max}_{qp}.
\]
\end{thm}

\begin{proof}
Notice that $\hat C=Q^TU^TCUQ$ is diagonal and hence
\begin{align*}
x(\xi)^TCx(\xi) & = \frac{1}{\displaystyle \max_{0\leq k\leq m} \xi^T\hat
A_k\xi} \xi^T Q^TU^TCUQ\xi = \frac{1}{\displaystyle \max_{0\leq k\leq m}
\xi^T\hat A_k\xi} \Tr(CX) .
\end{align*}
Therefore for any $\af >1 $ we have
%\begin{align*}
\begin{eqnarray*}
& & \prob \left\{ x(\xi)^T C x(\xi) \geq \frac{1}{\af}\, \Tr(CX) \right\} \\
& = & \prob \left\{ \max_{0\leq k \leq m} \xi^T\hat A_k\xi \leq \af  \right\} \\
& = & 1 - \prob \left\{ \max_{0\leq k \leq m} \xi^T\hat A_k\xi > \af  \right\} \\
& \geq & 1 - \prob \left\{ \max_{1\leq k\leq m} \xi^T\hat A_k \xi
> \af\right\} - \prob\left\{ \xi^T\hat A_0 \xi > \af \right\}.
%\end{align*}
\end{eqnarray*}
Since $\Tr(A_0) \leq 1$ and so $\af - \Tr(A_0) \geq 0$, it follows from
Lemma~\ref{qformprob} that
\[
\prob \left\{ \xi^T\hat A_0 \xi > \af \right\} \leq \prob \left\{
\sum_{1\leq i < j \leq m} (\hat A_0)_{ij} \xi_i \xi_j > 0\right\}
 < 1-\frac{1}{87}.
\]
Since $\hat A_k\succeq 0$ for $k=1,...,m$, and $\Tr(\hat A_k) \leq
1$, it follows from (12) in~\cite{NRT} that
\[
\prob \left\{\max_{1\leq k \leq m} \xi^T \hat A_k \xi > \af
\right\} < 2m \mu  e^{-\half \af}.
\]
Hence we have
\[
\prob \left\{ x(\xi)^T C x(\xi) \geq \frac{1}{\af}\, \Tr(C X)
\right\}
> \frac{1}{87} - 2m \mu  e^{-\half \af}.
\]
Letting $\af = 2\log(174\,m\mu)$ ensures the above probability to
be positive. Therefore, there exists a random vector $\xi$ such
that $\Tr(CX) \leq  \af\, x(\xi)^T C x(\xi)$, and the theorem is
proven.
%Therefore the quality bound in the theorem holds.
\end{proof}

We point out that Theorem~\ref{max qp bound} is an improvement of the
so-called approximate $S$-Lemma of
Ben-Tal, Nemirovski, and Roos~\cite{BNR02} (Lemma~A.6). In particular,
Ben-Tal~\emph{et al.}\ showed that $\alpha\le 2\log(16n^2\,m\mu)$, in contrast
to our bound $\alpha\le 2\log(174\,m\mu)$.

Notice that in (\ref{Qp-max}) there is only one indefinite inequality. A
natural question arises: can we allow more? %If one desires the
%bound on the ratio to be independent of the problem data (i.e., it
%depends only on the problem dimensions), then such ratio would be
%impossible in general when there are more than two indefinite
%inequality constraints, as the following example shows.
The following example shows that the answer is ``no" if
we wish to have a data-independent worst-case approximation ratio.
(Data-dependent approximation ratio bounds will be discussed
in Section~\ref{GHQP} where we do allow multiple indefinite constraints.)

\begin{example} \label{QP-max-3constraints}
Consider
\[
\begin{array}{rl} \max & x_1^2 + \frac{1}{M} x_2^2 \\
\st &   M x_1 x_2 + x_2^2 \le 1 \\
    & - M x_1 x_2 + x_2^2 \le 1 \\
    & M (x_1^2 - x_2^2) \le 1 ,
\end{array}
\]
where $M>0$ is an arbitrarily large positive constant. Its SDP
relaxation is
\[
\begin{array}{rl} \max & X_{11} + \frac{1}{M} X_{22} \\
\st & M X_{12} + X_{22} \le 1 ,\, - M X_{12} + X_{22} \le 1 ,\, M (X_{11} - X_{22}) \le 1 \\
& \left[ \begin{array}{cc} X_{11} & X_{12} \\ X_{21} & X_{22}
\end{array} \right] \succeq 0.
\end{array}
\]
For this quadratic program, the first two constraints imply that
$|x_1x_2| \le \frac{1- x_2^2}{M} \le \frac{1}{M}$ and so $x_1^2
\le \frac{1}{M^2 x_2^2}$. The third inequality assures that $x_1^2
\le \frac{1}{M}+x_2^2$. Therefore, $x_1^2 \le \min
\left\{\frac{1}{M^2 x_2^2},\frac{1}{M}+x_2^2\right\}\le
\frac{\sqrt{5}+1}{2M} \approx \frac{1.618}{M}$. Moreover, $x_2^2
\le 1$, and so $v^{\max}_{qp} \le \frac{2.618}{M}$.

The SDP relaxation satisfies both primal and dual Slater
conditions, so the primal-dual optimal solutions exist.
A feasible solution for the SDP relaxation (primal
problem) is the 2 by 2 identity matrix, with the objective value
being $1+\frac{1}{M}>1$. On the other hand, since $X_{22} \le
M|X_{12}|+X_{22} \le 1$, and $X_{11} \le X_{22} + \frac{1}{M}$, an
upper bound for the SDP optimal value is $1+\frac{2}{M}$.
Therefore, for this example, the ratio
$\frac{v^{\max}_{sdp}}{v^{\max}_{qp}} \ge \frac{M}{2.618} \approx
0.382 M$, which can be arbitrarily large, depending on the size of
$M$.
\end{example}

If there are at most two homogeneous quadratic constraints, and
moreover if the SDP relaxation has a primal-dual complementary
optimal solution, then the SDP optimal value will be equal to the
optimal value of the quadratic model; see e.g.\ Ye and
Zhang~\cite{YZ03} (Corollary~2.6). In other words, if there are no
more than two inequality constraints, then under the primal-dual
Slater condition, we will have $v^{\max}_{sdp}/v^{\max}_{qp}=1$.
In this sense, Example~\ref{QP-max-3constraints} is the smallest
possible in size. By removing the requirement that the SDP
relaxation has a finite optimal value, then it is possible to
construct an example which involves only two inequality
constraints.

\begin{example} \label{QP-max-2constraints}
Consider
\[
\begin{array}{rl} \max & x_1 x_2 + x_1^2 \\
\st &   x_1 x_2 \le 1 \\
    &   x_1^2 - x_2^2 \le 1,
\end{array}
\]
with the SDP relaxation
\[
\begin{array}{rl} \max & X_{12} + X_{11} \\
\st &   X_{12} \le 1,\, X_{11} - X_{22} \le 1, \\
 & \left[ \begin{array}{cc} X_{11} & X_{12} \\ X_{21} & X_{22}
\end{array} \right] \succeq 0.
\end{array}
\]

In terms of polar coordinates, $(x_1,x_2) \longrightarrow (r
\cos\theta, r \sin\theta)$, the original quadratic problem can be
turned into
\[
\begin{array}{rl} \max & r^2 ( \sin 2 \theta + \cos 2 \theta + 1)/2 \\
\st &   r^2 \sin 2 \theta  \le 2 \\
    &   r^2 \cos 2 \theta  \le 1.
\end{array}
\]
By a further change of variables
%$(r^2,2\theta) \longrightarrow (\rho, \alpha)$ and
$(r^2 \cos 2\theta, r^2 \sin 2\theta)\longrightarrow(y_1,y_2) $,
we can reformulate the original quadratic problem as
\[
\begin{array}{rl} \max & \frac{1}{2} \left( y_1 + y_2 + \sqrt{y_1^2 + y_2^2} \right) \\
\st &   y_1 \le 2 \\
    &   y_2 \le 1.
\end{array}
\]
This optimization problem has a unique optimal solution at $(y_1^*,y_2^*)=(2,1)$
with the optimal value being $\frac{3+\sqrt{5}}{2} \approx 2.618$.
The SDP relaxation problem is clearly unbounded, as any positive
multiple of the identity matrix is feasible. Therefore,
$v^{\max}_{sdp}/v^{\max}_{qp}=+\infty$. This example is possible
because the dual of the SDP relaxation problem is infeasible.
\end{example}

\section{Quadratic Optimization with Multiple Indefinite Constraints} \label{GHQP}

\setcounter{equation}{0}

Unlike the minimization form (\ref{minqp}) for which the SDP approximation
ratio can be infinite when there are more than one indefinite constraints (see
Example~\ref{QP min 4 constraints}), the maximization form (\ref{maxqp}) can
still admit a finite SDP approximation ratio in this case. In particular,
consider a general homogeneous quadratic maximization problem
\begin{equation} \label{max qp ind}
\begin{array}{rl}
\max & x^T C x \\
\st  & x^T A_k x \leq 1,\, k=0,1,...,m \\
     & x \in \F^n.
\end{array}
\end{equation}
Suppose that $\II,\DD$ are two index sets, $\II\cup
\DD=\{0,1,...,m\}$ and $\II\cap \DD=\emptyset$, such that $A_k
\succeq 0$ for $k\in \DD$ and $A_k$ indefinite for $k\in \II$. The
SDP relaxation for \reff{max qp ind} is
\begin{equation} \label{max qp sdp relaxation}
\begin{array}{rl}
\max & \Tr (C X)  \\
\st  & \Tr (A_k X) \leq 1,\, k=0,1,...,m \\
     & X \succeq 0.
\end{array}
\end{equation}

We begin our analysis with a technical lemma which bounds the
probability of an exponential tail. Similar bounds exist in the
literature, e.g.~\cite{BobKov}. However, the lemma below serves
our needs exactly; for completeness we include a proof here.

\begin{lemma} \label{chi2 posi}
Let $\{\lambda_i\}_{i=1}^{n}$ be any given real numbers and
 $\{\eta_i\}_{i=1}^n$ be i.i.d.\
random variables drawn from either the real or complex valued zero mean Gaussian
distribution with unit variance. Let $\sigma = \sqrt{\sum_{i=1}^n \lambda_i^2}$ and
$\delta=\max\left\{ \max\{\lambda_i\mid 1\le i\le n\},0\right\}$.
Then, for any $\alpha>0$ there holds
\[
\prob \left\{ \sum_{i=1}^n \lambda_i \eta_i^2 - \sum_{i=1}^n
\lambda_i \ge \alpha \sigma \right\} \le \left\{
\begin{array}{ll}
\exp
\left(-\min\left\{\alpha,\frac{\sigma}{\delta}\right\}\frac{\alpha}{8}\right),
& \mbox{if $\eta_i\sim N(0,1)$ is real Gaussian,}\\ [10pt]
\exp
\left(-\min\left\{\alpha,\frac{\sigma}{\delta}\right\}\frac{\alpha}{4}\right),
&\mbox{if $\eta_i\sim N_c(0,1)$ is complex Gaussian.}
\end{array}
\right.\]
%
%Suppose that $\lambda_i \in \re$, $i=1,...,n$, and $\eta_i^2$ are
%i.i.d exponential distributions with mean $1$. Let $\sigma =
%\sqrt{\sum_{i=1}^n \lambda_i^2}$ and $\delta=\max\left\{
%\max\{\lambda_i\mid 1\le i\le n\},0\right\}$. Then, for any
%$\alpha>0$ there holds
%\[
%\prob \left\{ \sum_{i=1}^n \lambda_i \eta_i^2 - \sum_{i=1}^n
%\lambda_i \ge \alpha \sigma \right\} \le \exp
%\left(-\min\left\{\alpha,\frac{\sigma}{\delta}\right\}\frac{\alpha}{4}\right).
%\]
\end{lemma}
\begin{proof}
We will only prove the real Gaussian case; the complex case is similar and therefore
omitted. Let
$\beta:=\frac{1}{4}\min\{\frac{1}{\delta},\frac{\alpha}{\sigma}\}$.
Then, $2\beta \lambda_i \le 1/2$ for all $i=1,...,n$, and $\beta
\sigma =\frac{1}{4}\min\{\frac{\sigma}{\delta},\alpha\}$. Note
that for any $t \le 1/2$ the following inequality holds:
\begin{equation} \label{simple ineq}
\frac{1}{1-t} \le e^{t+t^2}.
\end{equation}
Let $ \zeta:= e^{\beta \sum_{i=1}^n \lambda_i \eta_i^2}$. Since
$\{\eta_i^2\}_{i=1}^n$ are standard i.i.d.\ $\chi^2$ random
variables, it follows that
\[
%\begin{eqnarray*}
\ex (\zeta) = \prod_{i=1}^n \ex \left(e^{\beta \lambda_i
\eta_i^2}\right) = \prod_{i=1}^n \frac{1}{\sqrt{1-2\beta
\lambda_i}} = \left( \prod_{i=1}^n \frac{1}{1-2\beta \lambda_i}
\right)^{\frac{1}{2}} \le \left( \prod_{i=1}^n e^{2 \beta
\lambda_i + 4 \beta^2 \lambda_i^2} \right)^{\frac{1}{2}}   =
e^{  2 \beta^2 \sigma^2+ \beta \sum_{i=1}^n \lambda_i   }
%\end{eqnarray*}
\]
where the inequality is due to \reff{simple ineq}. This together with the
Markov inequality implies
\begin{eqnarray*}
\prob \left\{ \sum_{i=1}^n \lambda_i \eta_i^2 - \sum_{i=1}^n
\lambda_i \ge \alpha \sigma \right\} &=& \prob \left\{ \zeta \ge
e^{  \beta (  \alpha \sigma+\sum_{i=1}^n \lambda_i )} \right\} \\
&\le& \frac{\ex (\zeta)}{e^{  \beta (  \alpha \sigma+\sum_{i=1}^n \lambda_i )} } \\
&\le& e^{2\beta^2 \sigma^2 - \beta \sigma \alpha}  = e^{\beta
\sigma( 2 \beta \sigma - \alpha)} \le
%e^{-\min\left\{\alpha,\frac{\sigma}{\delta}\right\}\frac{3\alpha}{8}}
e^{\beta \sigma( \frac{\alpha}{2} - \alpha)} \\
&=&
e^{-\min\left\{\alpha,\frac{\sigma}{\delta}\right\}\frac{\alpha}{8}}.
\end{eqnarray*}
The lemma is proven.
\end{proof}

We are now ready to pursue the performance analysis for the real case $\F=\re$.
Assume that \reff{max qp sdp relaxation} has an optimal solution
$\hat X$. Denote the set of (real) eigenvalues of $A_k \hat X$ as
$\lambda^k_1,...,\lambda^k_n$, $k=0,1,...,m$. Since $\Tr (A_k \hat
X) \le 1$, it follows that $\sum_{i=1}^n \lambda^k_i \le 1$.
Moreover, $\|A_k \hat X\|_F^2 = \sum_{i=1}^n (\lambda^k_i)^2$,
$k=0,1,...,m$, where $\|\cdot\|_F$ denotes the Frobenius norm of a
matrix.
%For a given
%$\alpha>1$, we consider the probability of the joint event
%%\begin{equation} \label{joint event}
%\[
%\left\{ \begin{array}{l} \xi^T A_k \xi \le \alpha,\, k=1,...,m \\
%\xi^T C \xi \ge \Tr C \hat X . \end{array} \right.
%\]
%%\end{equation}

Let $\xi$ be a random
vector drawn from the Gaussian distribution $N(0,\hat X)$.
For any $\alpha> 1$ and $0\le k\le m$, we
consider the probability of the event
$\prob \{ \xi^T A_k \xi
> \alpha\}$. By diagonalization, we have $\prob \{ \xi^T
A_k \xi > \alpha\} = \prob \{ \sum_{i=1}^n \lambda^k_i \eta_i^2
> \alpha\} $, where $\eta=(\eta_1,...,\eta_n)^T$ is a random
vector following the normal distribution $N(0,I_n)$.

If we let $\sigma^k:=\sqrt{\sum_{i=1}^n
(\lambda^k_i)^2}=\| A_k \hat X\|_F$, and
$\delta^k:=\max\left\{0,\max\{\lambda^k_i\mid 1\le i\le
n\}\right\}$, then Lemma~\ref{chi2 posi} leads to
\begin{equation} \label{psd part}
\prob \{ \xi^T A_k \xi
> \alpha\} \le \exp \left(-\min\left\{\frac{\alpha-\sum_{i=1}^n
\lambda^k_i}{\sigma^k},\frac{\sigma^k}{\delta^k}\right\}\frac{\alpha-\sum_{i=1}^n
\lambda^k_i}{8 \sigma^k}\right),\quad \forall\; 0\le k\le m.
\end{equation}
%for any $0\le k\le m$.
Alternatively, we can bound the tail probability using Chebyshev's inequality.
In particular, since
$\Var(\sum_{i=1}^n \lambda^k_i \eta_i^2)=2 \sum_{i=1}^n
(\lambda^k_i)^2 = 2 \| A_k \hat X\|_F^2$, it follows from Chebyshev's
inequality
\begin{eqnarray}
\prob \left\{ \sum_{i=1}^n \lambda^k_i \eta_i^2
> \alpha \right\} &=& \prob \left\{ \sum_{i=1}^n \lambda^k_i \eta_i^2
- \sum_{i=1}^n \lambda^k_i > \alpha - \sum_{i=1}^n \lambda^k_i
\right\} \nonumber \\
&\le& \prob \left\{ \left|\sum_{i=1}^n \lambda^k_i \eta_i^2 -
\sum_{i=1}^n \lambda^k_i \right| > \alpha - \sum_{i=1}^n
\lambda^k_i
\right\}  \nonumber \\
&\le & \frac{\Var( \sum_{i=1}^n \lambda^k_i \eta_i^2
)}{\left(\alpha - \sum_{i=1}^n \lambda^k_i\right)^2} %\nonumber \\
\le \frac{2 \| A_k \hat X\|_F^2}{(\alpha-1)^2},\quad \forall\; 0\le k\le m, \label{prob est}
\end{eqnarray}
where we have used the fact $\alpha>1\ge \sum_{i=1}^n \lambda^k_i$.
Applying Lemma~\ref{realpb} and using
\reff{prob est}--\reff{psd part} gives
\begin{eqnarray*}
& & \prob \left\{ \xi^T A_k \xi \le \alpha,\, k=0,1,...,m;\, \xi^T C
\xi \ge \Tr (C \hat X) \right\} \\
&\ge& 1 - \prob \left\{\xi^T C \xi < \Tr (C \hat X) \right\} -
\sum_{k=0}^m\prob \left\{ \xi^T A_k \xi > \alpha\right\} \\
&\ge& \frac{3}{100} - \sum_{k=0}^m \min\left\{
\exp\left(-\min\left\{\frac{\alpha-\sum_{i=1}^n
\lambda^k_i}{\sigma^k},\frac{\sigma^k}{\delta^k}\right\}\frac{\alpha-\sum_{i=1}^n
\lambda^k_i}{8\sigma^k}\right), \frac{2 \| A_k \hat
X\|_F^2}{(\alpha-1)^2}\right\} .
\end{eqnarray*}
Notice that $\delta^k\le \sigma^k$ and $\sum_{i=1}^n \lambda^k_i\le 1$
for any $k$. Therefore, we have, for any $\alpha>1$,
\begin{eqnarray*}
& & \prob \left\{ \xi^T A_k \xi \le \alpha,\, k=0,1,...,m;\, \xi^T C
\xi \ge \Tr (C \hat X) \right\} \\
&\ge& \frac{3}{100} - \sum_{i \in \DD}
\exp\left(-\min\left\{\frac{\alpha-1}{\sigma^k},1\right\}\frac{\alpha-1}{8\sigma^k}\right) \\
& & -\sum_{i\in \II} \min\left\{
\exp\left(-\min\left\{\frac{\alpha-1}{\sigma^k},1\right\}\frac{\alpha-1}{8\sigma^k}\right),
\frac{2 \| A_k \hat X\|_F^2}{(\alpha-1)^2}\right\}  .
\end{eqnarray*}

Let us choose
\[
\alpha=1+\max\left\{ 20 + 8 \log |\DD|, \min\left\{ (20 + 8 \log
|\II|)\max_{k \in \II} \|A_k \hat X\|_F,\sqrt{200\sum_{k\in \II}
\|A_k \hat X\|_F^2} \right\}\right\}.
%\alpha\ge 1+ (20 + 8 \log |\DD|)
\]
Since $\sigma^k \le \sum_{i=1}^n \lambda^k_i\le 1$ for $k\in \DD$, it follows from
the choice of $\alpha$ that
\[
\exp\left(-\min\left\{\frac{\alpha-1}{\sigma^k},1
\right\}\frac{\alpha-1}{8\sigma^k}\right) =
\exp\left(-\frac{\alpha-1}{8\sigma^k}\right) \le
\exp\left(-\frac{\alpha-1}{8}\right) \le \frac{1}{100 |\DD|}, \quad\forall\; k\in \DD,
\]
%for any $k\in \DD$, and
and
\[
\sum_{i\in \II} \min\left\{ \exp
\left(-\min\left\{\frac{\alpha-1}{\sigma^k},1\right\}\frac{\alpha-1}{8\sigma^k}\right),
\frac{2 \| A_k \hat X\|_F^2}{(\alpha-1)^2}\right\}\le
\frac{1}{100}.
\]
This further implies that
\[
\prob \left\{ \xi^T A_k \xi \le \alpha,\, k=0,1,...,m;\, \xi^T C \xi
\ge \Tr (C \hat X) \right\}\ge \frac{1}{100}.
\]

Summarizing, we obtain the following worst-case performance ratio
bounds on the SDP relaxation for a real-valued homogeneous
(indefinite) quadratic maximization problem. [We also state the complex case without proof.]

\begin{theorem} \label{max qp sdp ind}
For the quadratic optimization problem $\reff{max qp ind}$ with
$\F=\re$ and its SDP relaxation $\reff{max qp sdp relaxation}$,
suppose that an optimal solution, say $\hat X$, for $\reff{max qp sdp relaxation}$
exists. Then,
\[
\frac{v^{\max}_{sdp}}{v^{\max}_{qp} } \le 1+\max\left\{ 20+8 \log
|\DD| , \min\left\{(20+8 \log |\II|)\max_{k \in \II} \|A_k \hat
X\|_F,\sqrt{200\sum_{k\in \II} \|A_k \hat X\|_F^2}
\right\}\right\}.
\]
Similarly, for the complex case $\F=\cpx$, % and its SDP relaxation $\reff{max qp sdp
%relaxation}$, %suppose that an optimal solution for $\reff{max qp sdp
%relaxation}$ exists and let it be $\hat X$. Then,
we have
\[
\frac{v^{\max}_{sdp}}{v^{\max}_{qp} } \le 1+\max\left\{15+4 \log
|\DD| , \min\left\{(15+4 \log |\II|)\max_{k \in \II} \|A_k \hat
X\|_F,\sqrt{40\sum_{k\in \II} \|A_k \hat X\|_F^2}
\right\}\right\}.
\]
\end{theorem}

Let us consider two special cases of Theorem~\ref{max qp sdp ind}.
First, if $\II=\emptyset$, then Theorem~\ref{max qp sdp ind}
becomes $\frac{v^{\max}_{sdp}}{v^{\max}_{qp} } \le 20+8 \log m $
(in the real case), which recovers the approximation result   of
Nemirovski
{\em et al.}~\cite{NRT}. %, and Theorem~\ref{min qp sdp ind}
%recovers an approximation result in Luo {\em et al.}~\cite{LSTZ}.
The second case is $\DD=\emptyset$, where Theorem~\ref{max qp sdp ind} becomes
$$\frac{v^{\max}_{sdp}}{v^{\max}_{qp} } \le 1+ \min\left\{(20+8 \log(m+1))\max_{0\le k\le m} \|A_k \hat
X\|_F,\sqrt{200\sum_{k=0}^m \|A_k \hat X\|_F^2}
\right\}.$$ Below is an example showing that
this bound is also
tight (in the order of magnitude). Specifically, consider
Example~\ref{QP-max-3constraints} again:
\[
\begin{array}{rl} \max & x_1^2 + \frac{1}{M} x_2^2   \\
\st &   M x_1 x_2 + x_2^2 \le 1 \\
    & - M x_1 x_2 + x_2^2 \le 1 \\
    & M (x_1^2 - x_2^2) \le 1 .
\end{array}
\]
In this case we know that the SDP relaxation has an optimal
solution $\hat X=\left[ \begin{array}{cc} 1+ \frac{1}{M} & 0 \\
0 & 1\end{array}\right]$, while the approximation ratio is
$v^{\max}_{sdp}/v^{\max}_{qp} = O(M)$. There are three
constraints, all indefinite, $\II=\{1,2,3\}$, with
\[
A_1=\left[ \begin{array}{cc} 0 & \frac{M}{2} \\ \frac{M}{2} & 1
\end{array}\right],\,
A_2=\left[ \begin{array}{cc} 0 & -\frac{M}{2} \\ -\frac{M}{2} & 1
\end{array}\right],\,
A_3=\left[ \begin{array}{cc} M & 0 \\ 0 & -M
\end{array}\right],
\]
and so one may compute that
\[
A_1\hat X=\left[ \begin{array}{cc} 0 & \frac{M}{2} \\
\frac{M}{2}+\frac{1}{2} & 1
\end{array}\right],\,
A_2\hat X=\left[ \begin{array}{cc} 0 & -\frac{M}{2} \\
-\frac{M}{2} -\frac{1}{2} & 1
\end{array}\right],\,
A_3 \hat X=\left[ \begin{array}{cc} M+1 & 0 \\ 0 & -M
\end{array}\right].
\]
Thus, $\| A_k \hat X\|_F^2=O(M^2)$, for $k=1,2,3$.
Theorem~\ref{max qp sdp ind} predicts that
$v^{\max}_{sdp}/v^{\max}_{qp} \le O(M)$, and this upper bound is
exactly attained
in this example. %Therefore, the order of magnitude of the bound
%presented in Theorem~\ref{max qp sdp ind} cannot be improved in
%general.

\section{Simulations and Discussions} \label{conclusion}

This paper studies the quality bound of SDP relaxation for solving
nonconvex quadratic optimization problems~\reff{hqp1}
and~\reff{Qp-max}. For problem~\reff{hqp1}, a quality
bound $O(m^2)$ was derived for $\F=\re$, and a quality bound
$O(m)$ for $\F=\cpx$, when there is only one constraint
$x^*A_0x\geq 1$ with $A_0$ indefinite. For problem~\reff{Qp-max},
a quality bound $O(\log m)$ was derived when there is only one
nonconvex constraint $x^*A_0x\leq 1$ with $A_0$ indefinite. These
quality bounds are independent of the problem dimension $n$ or data matrices,
and only depend on the number of constraints.

For problem~\reff{hqp1}, if there are two or more constraints in
the form of $x^*A x\geq 1$ with $A$ indefinite, then there is no
general quality bound as shown by Example~\ref{QP min 4
constraints}. For problem~\reff{Qp-max}, if there are two or more
nonconvex constraints, a quality bound is given in
Theorem~\ref{max qp sdp ind}, albeit the bound is dependent not
only on the number of constraints but also the data of the
problem.

As shown in the preceding sections, these quality bounds are
derived based on the worst-case analysis, and they are indeed
tight, in the worst case, up to some constant. This analysis is
important as a theoretical guide. The empirical tests, on the
other side, serve a quite different purpose. Next, we present some
numerical experiments on randomly
generated instances. These numerical experiments show
that the average approximation ratios are
much better than the worst-case ratio, even though they appear to still
follow the same  growth trend (as a function of $m$).

\begin{figure}[tb]
%\btab{@{}c@{}c@{}}
\hspace*{1cm}\includegraphics[height=.5\textheight]{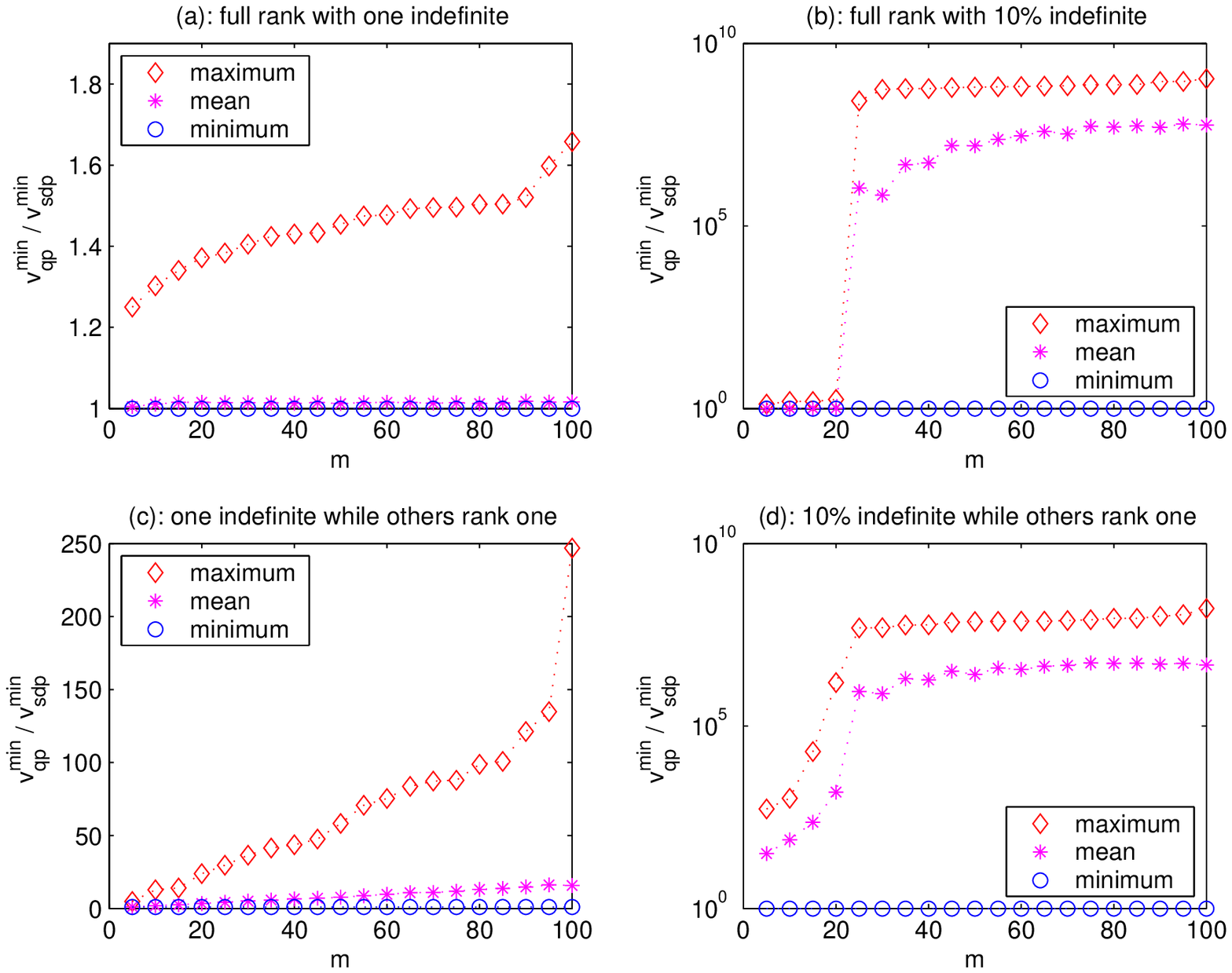}
%\etab
\caption{Empirical quality bounds for problem \reff{hqp1} with $\F=\re$ } \label{ratmin}
\end{figure}

%In this section, we use random examples to test the performance of
%SDP relaxations. V
More specifically, we generate various random symmetric matrices $A_k$
 in the following way: for a full rank positive
semidefinite $A_k$, we set $A_k = \mbox{rand} \cdot Q^T\cdot
\mbox{diag}(\mbox{abs}(\mbox{randn}(n,1)))\cdot Q$, where
`$\mbox{rand}$', `$\mbox{randn}$' are Matlab notations,
and $Q$ is an orthogonal matrix obtained by QR factorization of a
random matrix $\mbox{randn}(n)$; for a rank-one positive
semidefinite $A_k$, we set $A_k = \mbox{rand} \cdot Q^T \cdot
\mbox{diag}([\mbox{abs}(\mbox{randn}); \mbox{zeros}(n-1,1)]))\cdot
Q$; and for an indefinite
$A_k$, we set $A_k = \mbox{rand} \cdot Q^T \cdot
\mbox{diag}(\mbox{randn}(n,1))\cdot Q$ ($Q$ defined as before).
To examine the performance of SDP relaxation for randomly
generated problem of form \reff{hqp1} with $\F=\re$,
%Since the quality bounds are independent of problem size $n$, fix
 we simply set $C$ to be the identity matrix: fix
$n=10$, and choose $m$ from $5,10,15,...,100$. For each $m$, we do
the following: (a) Generate $1,000$ random problems such that only
one of the $A_k$'s is indefinite, and all the other $A_k$'s are
positive definite; (b) Generate $1,000$ random problems such that
$10\%$ of the $A_k$'s are indefinite, and all the other $A_k$'s
are positive definite; (c) Generate $1,000$ random problems such
that only one of the $A_k$'s is indefinite, and all the other
$A_k$ are rank one and positive semidefinite; (d) Generate $1,000$
random problems such that $10\%$ of the $A_k$'s are indefinite,
and all the other $A_k$ are rank one and positive semidefinite.
For each instance of the above randomly generated problems, we
solve its SDP relaxation to obtain an optimal solution $Z^*$ and
optimal value $v_{sdp}^{\min}$. Then we find one approximate
solution for \reff{hqp1} by the following randomization process.
Generate $100$ random vectors $\xi^1,...,\xi^{100}$. For each
$k=1,...,100$, let $x^i =
\xi^i/\sqrt{\min_{k=1}^m(\xi^i)^TA_k\xi^i}$. Then
$v_{qp}^{\min}\leq \hat v_{qp}^{\min}:=\min_{i=1}^{100}
(x^i)^TCx^i$. We use the empirical quality bound $\hat
v_{qp}^{\min}/ v_{sdp}^{\min}$ to estimate the real quality bound
$v_{qp}^{\min}/v_{sdp}^{\min}$, since the former is at least the
latter. These empirical quality bounds are plotted in
Figure~\ref{ratmin}. \bcen
\ecen
In this figure, diamonds
$\diamond$ are the maximum quality bounds in the $1,000$ random
problems for each $m$, stars $*$ are the mean quality bounds, and
circles $\circ$ are the minimum quality bounds. %The upper left
%plot is for case (a), the upper right plot for case (b), the lower
%left plot for case (c), and the lower right plot for case (d).
For cases (a) and (c), we have quality bound $\mathcal{O}(m^2)$
for the worst case, while for cases (b) and (d), there is no
worst-case theoretical quality bound. As we can see, the computed
empirical bounds are very small for case (a), and are moderate for
case (c), and are indeed big for cases (b) and (d).
%This is consistent with our theoretical analysis.
%For the case of complex field $\F=\cpx$, we have similar empirical numerical performance.

\begin{figure}[htb]
%\btab{@{}c@{}c@{}}
\hspace{1cm}
\includegraphics[height=.5\textheight]{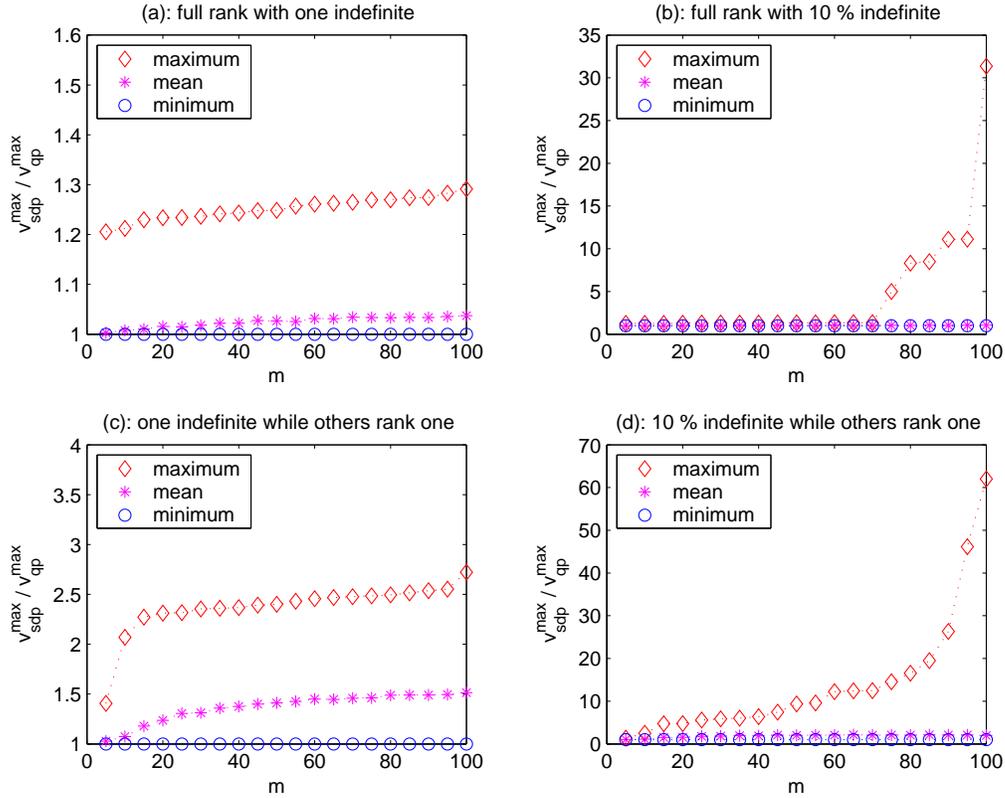}
%\etab
\caption{Empirical quality bounds for problem \reff{Qp-max}
} \label{ratmax}
\end{figure}

To examine the SDP relaxation performance for the maximization
problem~\reff{Qp-max}, we generate four classes of random problems
in the same way as for problem~\reff{hqp1}, except that the matrix
$C$ in the objective is now indefinite (generated the same way as
an indefinite $A_k$). After the SDP relaxation is solved, we apply
a similar randomization procedure to find a lower bound $\hat
v_{qp}^{\max}$ for $v_{qp}^{\max}$ as we did for
problem~\reff{hqp1}. The empirical quality bound
$v_{sdp}^{\max}/\hat v_{qp}^{\max}$ is an upper bound for the
actual quality bound $v_{sdp}^{\max}/v_{qp}^{\max}$. The empirical
quality bounds are plotted in Figure~\ref{ratmax}. The legends
$\diamond$, $*$ and $\circ$ carry the same meaning as they did in
Figure~\ref{ratmin}.
%The upper left plot is for case (a), the
%upper right plot for case (b), the lower left plot for case (c),
%and the lower right plot for case (d).
This figure shows that the computed empirical bounds are close to
one for cases (a) and (c), and are somewhat larger for cases (b)
and (d). This is consistent with the bounds in Theorem~\ref{max qp
sdp ind}.

\bigskip
\noindent {\bf Acknowledgement:} The authors wish to thank Yuval
Peres for suggesting the reference~\cite{KS} to us.

\bcen
\ecen

\end{document}